\documentclass[12pt]{article}
\usepackage{amsfonts}
\usepackage{mathrsfs}

\usepackage{amsmath,amssymb,amsfonts,theorem,makeidx,latexsym,epsfig,color,subfigure,algorithm,algorithmic}\usepackage[numbers,sort&compress]{natbib}

\newtheorem{defn}{Definition}[section]

\newtheorem{lem}[defn]{Lemma}

{\theorembodyfont{\rmfamily}

}

\newtheorem{thm}[defn]{Theorem}

\newtheorem{obs}[defn]{Observation}

\newtheorem{cor}[defn]{Corollary}

\numberwithin{equation}{section}

\def\bp{{\noindent\bf Proof. \ }}

\def\ep{\hfill$\square$\par\bigskip}

\title{Disjunctive domination in trees \footnote{The research is supported by NSFC (No. 11301440),
Natural Science Foundation of Fujian Province (CN)(2015J05017)}}

\author{{Wei Zhuang$^{a}$\thanks{Corresponding author; E-mail: zhuangweixmu@163.com}\quad Litao Guo$^a$ \quad Guoliang Hao $^b$\vspace*{0.8cm}}\\
        {\small $^a$ School of Applied Mathematics, Xiamen University of Technology,}\\
        {\small Xiamen Fujian 361024, P.R.China \vspace*{0.3cm}}\\
        {\small $^b$ College of Science, East China University of Technology,}\\
        {\small Nanchang Jiangxi 330013, P.R.China}}

\date{}

\begin{document}

\maketitle

\begin{abstract} In this paper, we study a
parameter that is a relaxation of arguably the most important
domination parameter, namely the domination number. Given the sheer
scale of modern networks, many existing domination type structures
are expensive to implement. Variations on the theme of dominating
sets studied to date tend to focus on adding restrictions which in
turn raises their implementation costs. As an alternative route a
relaxation of the domination number, called disjunctive domination,
was proposed and studied by Goddard et al. A set $D$ of vertices in
$G$ is a disjunctive dominating set in $G$ if every vertex not in
$D$ is adjacent to a vertex of $D$ or has at least two vertices in
$D$ at distance $2$ from it in $G$. The disjunctive domination
number, $\gamma^{d}_2(G)$, of $G$ is the minimum cardinality of a
disjunctive dominating set in $G$. We show that if $T$ is a tree of
order $n$ with $l$ leaves and $s$ support vertices, then
$\frac{n-l+3}{4}\leq \gamma^{d}_2(T)\leq \frac{n+l+s}{4}$. Moreover,
we characterize the families of trees which attain these bounds.
\end{abstract}

\begin{minipage}{150mm}


\end{minipage}

\section{Introduction}
Over the last few decades, the scale of networks and the role of
graphs as models for networks has changed, and in practical terms,
many existing domination type structures are too expensive to
implement. The majority of domination-type variants studied to date
tend to focus on adding restrictions which in turn raises their
implementation costs. As a result the idea of relaxing conditions on
domination-type parameters is appealing. A relaxation of the
domination number, called disjunctive domination, was proposed and
studied in \cite{Goddard}. In this paper we continue the study of
disjunctive domination in graphs.

A \emph{dominating set} in a graph $G$ is a set $S$ of vertices of
$G$ such that every vertex in $V(G)\setminus S$ is adjacent to at
least one vertex in $S$. The domination number of $G$, denoted by
$\gamma(G)$, is the minimum cardinality of a dominating set of $G$.
A set $D$ of vertices in a graph $G$ is a \emph{disjunctive
dominating set}, abbreviated \emph{$2DD$-set}, in $G$ if every
vertex not in $D$ is adjacent to a vertex of $D$ or has at least two
vertices in $D$ at distance $2$ from it in $G$. We say a vertex $v$
in $G$ is \emph{$2D$-dominated}, by the set $D$, if $N[v]\cap D\neq
\emptyset$ or there exist at least two vertices in $D$ at distance
$2$ from $v$ in $G$. The \emph{disjunctive domination number} of
$G$, denoted by $\gamma^{d}_2(G)$, is the minimum cardinality of a
$2DD$-set in $G$. A disjunctive dominating set of $G$ of cardinality
$\gamma^{d}_2(G)$ is called a $\gamma^{d}_2(G)$-set. If the graph
$G$ is clear from the context, we simply write $\gamma^{d}_2$-set
rather than $\gamma^{d}_2(G)$-set.

Every dominating set is a $2DD$-set. The concept of disjunctive
domination in graphs has been studied in \cite{Goddard, Henning1,
Henning2, Henning3} and elsewhere.

Let $G=(V, E)$ be a graph with vertex set $V$ of order $n(G)=|V|$
and edge set $E$ of size $m(G)=|E|$, and let $v$ be a vertex in $V$.
The \emph{open neighborhood} of $v$ is $N(v)=\{u\in V|uv\in E\}$ and
the \emph{closed neighborhood} of $v$ is $N[v]=N(v)\cup \{v\}$. The
\emph{degree} of a vertex $v$ is $d(v)=|N(v)|$. For two vertices $u$
and $v$ in a connected graph $G$, the \emph{distance} $d(u, v)$
between $u$ and $v$ is the length of a shortest $(u, v)$-path in
$G$. The maximum distance among all pairs of vertices of $G$ is the
\emph{diameter} of a graph $G$ which is denoted by $diam(G)$. A
\emph{leaf} of $G$ is a vertex of degree $1$ and a \emph{support
vertex} of $G$ is a vertex adjacent to a leaf. Denote the sets of
leaves and support vertices of $G$ by $L(T)$ and $S(T)$,
respectively. Let $l(T)=|L(T)|$ and $s(T)=|S(T)|$. A \emph{double
star} is a tree that contains exactly two vertices that are not
leaves.

\section{Main results}
In this paper, we give a lower bound and an upper bound for the
disjunctive domination number of a tree in terms of its order, the
number of leaves and support vertices in the tree. Further, we
provide the constructive characterization of trees that achieve
equality in the two bounds. We state this formally as follows.

\begin{obs}\cite{Henning3}
If $T$ is a tree of order at least $3$, then we can choose a
$\gamma^{d}_2$-set of $T$ contains no leaf.
\end{obs}

\begin{cor}
Let $T$ be a tree of order at least $3$ and $D$ be a
$\gamma^{d}_2$-set of $T$ contains no leaf, if a support vertex has
degree two, then it belongs to $D$.
\end{cor}

By a weak partition of a set we mean a partition of the set in which
some of the subsets may be empty. For our purposes, we define a
\emph{labeling} of a tree $T$ as a weak partition $S=(S_A, S_B, S_C,
S_D)$ of $V(T)$ (This idea of labeling the vertices is introduced in
\cite{Dorfling}). We will refer to the pair $(T, S)$ as a
\emph{labeled tree}. The label or \emph{status} of a vertex $v$,
denoted sta$(v)$, is the letter $x\in \{A, B, C, D\}$ such that
$v\in S_x$. Next, we ready to give two families $\mathscr{T}_1$ and
$\mathscr{T}_2$, each member of which is obtained from the labeled
trees $(P_3, S')$ and $(P_4, S'')$ respectively by a series of
operations. Before this, we give two definitions. If a labeled tree
$(T, S)\in \mathscr{T}_2$, the path $P_4$ (which comes from the
labeled tree $(P_4, S'')$) is an induced path of $T$, and we call it
the \emph{basic path} of $T$. For a vertex $v\not \in S(T)$, which
has status $A$ and does not belong to the basic path, if there
exists a vertex $u$ such that $vv_1v_2u$ is an induced path of $T$
and sta$(v_1)=C$, sta$(v_2)=D$, sta$(u)=B$, we call $u$ a
\emph{corresponding vertex} of $v$. In addition, for a vertex $u$,
which has status $B$, if there exists a vertex $v$ such that
$vv_1v_2u$ is an induced path of $T$ and sta$(v)=A$, sta$(v_1)=C$,
sta$(v_2)=D$, we call $v$ a \emph{corresponding vertex} of $u$.

\begin{center}
  \includegraphics[width=3in]{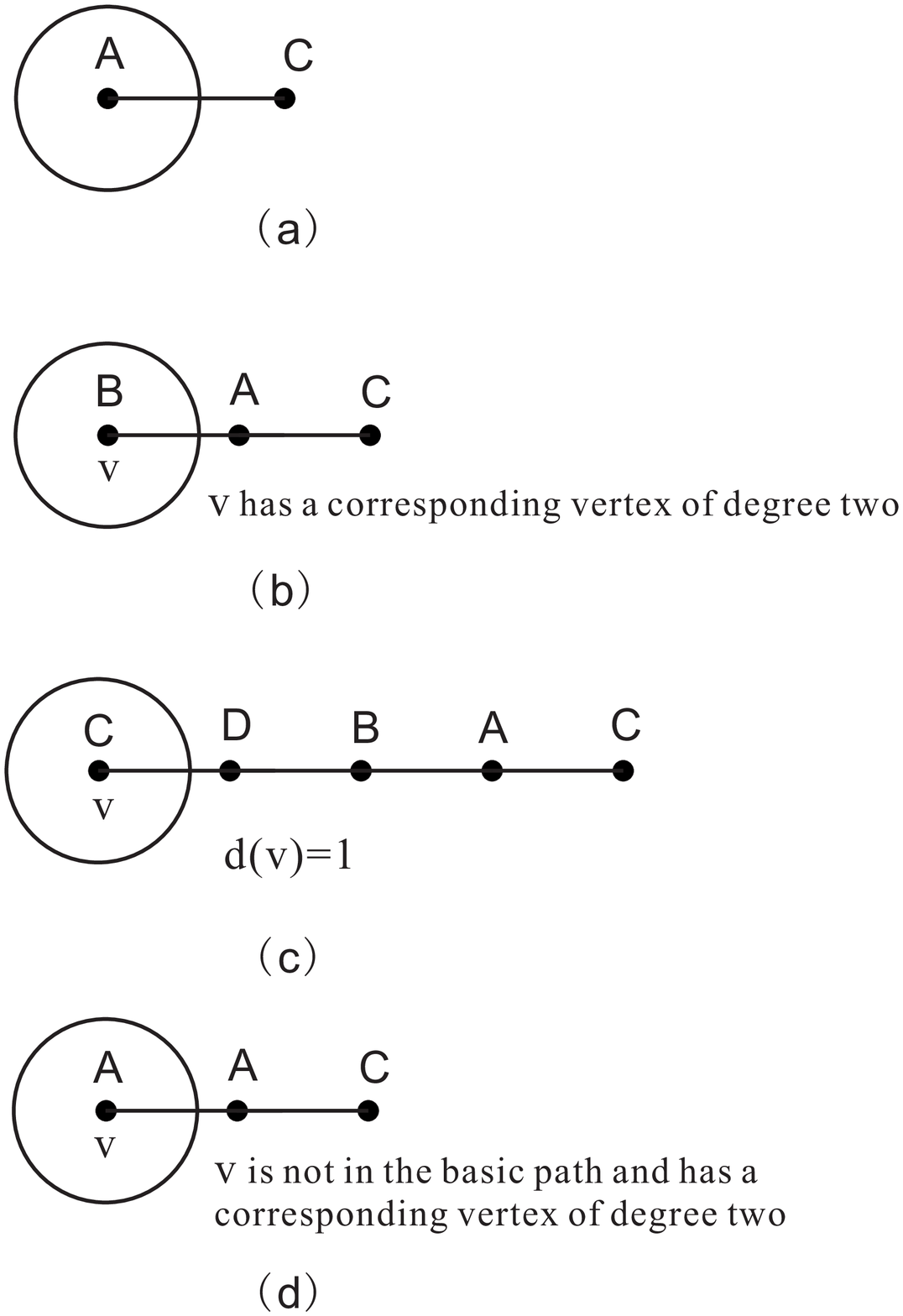}
 \end{center}
\qquad \qquad \qquad \qquad \qquad \qquad \qquad \qquad \qquad
{\small Fig.1}

In what follows, we give four operations as follows:

{\bf Operation} $\mathscr{O}_1$: Let $v$ be a vertex with
sta$(v)=A$. Add a vertex $u$ and the edge $uv$. Let sta$(u)=C$.

{\bf Operation} $\mathscr{O}_2$: Let $v$ be a vertex with sta$(v)=B$
that has a corresponding vertex of degree two. Add a path $u_1u_2$
and the edge $u_1v$. Let sta$(u_1)=A$, sta$(u_2)=C$.

{\bf Operation} $\mathscr{O}_3$: Let $v$ be a vertex with sta$(v)=C$
that has degree one. Add a path $u_1u_2u_3u_4$ and the edge $u_1v$.
Let sta$(u_1)=D$, sta$(u_2)=B$, sta$(u_3)=A$, sta$(u_4)=C$.

{\bf Operation} $\mathscr{O}_4$: Let $v$ be a vertex not in the
basic path that has status $A$ and has a corresponding vertex of
degree two. Add a path $u_1u_2$ and the edge $u_1v$. Let
sta$(u_1)=A$, sta$(u_2)=C$.

The three operations $\mathscr{O}_1$, $\mathscr{O}_2$,
$\mathscr{O}_3$ and $\mathscr{O}_4$ are illustrated in Fig.1(a),
(b), (c) and (d).

Let $\mathscr{T}_1$ be the minimum family of labeled trees that: (i)
contains $(P_3, S')$ and $S'$ is the labeling that assigns to the
two leaves of the path $P_3$ status $C$, and the central vertex
status $A$; and (ii) is closed under the two operations
$\mathscr{O}_1$ and $\mathscr{O}_3$ that are listed as above, which
extend the tree $T'$ to a tree $T$ by attaching a tree to the vertex
$v\in V(T')$.

Let $\mathscr{T}_2$ be the minimum family of labeled trees that: (i)
contains $(P_4, S'')$ where $S''$ is the labeling that assigns to
the two leaves of the path $P_4$ status $C$, and the remaining
vertices status $A$; and (ii) is closed under the three operations
$\mathscr{O}_2$, $\mathscr{O}_3$ and $\mathscr{O}_4$ that are listed
as above, which extend the tree $T'$ to a tree $T$ by attaching a
tree to the vertex $v\in V(T')$.

We take an example to make it easier for reader to understand the
family $\mathscr{T}_1$ and $\mathscr{T}_2$. The trees are depicted
in Fig.2(a) and (b) belong to $\mathscr{T}_1$ and $\mathscr{T}_2$,
respectively. In Fig.2(b), the induced path $v_1v_2v_3v_4$ is the
basic path of the tree.

\begin{center}
  \includegraphics[width=4.5in]{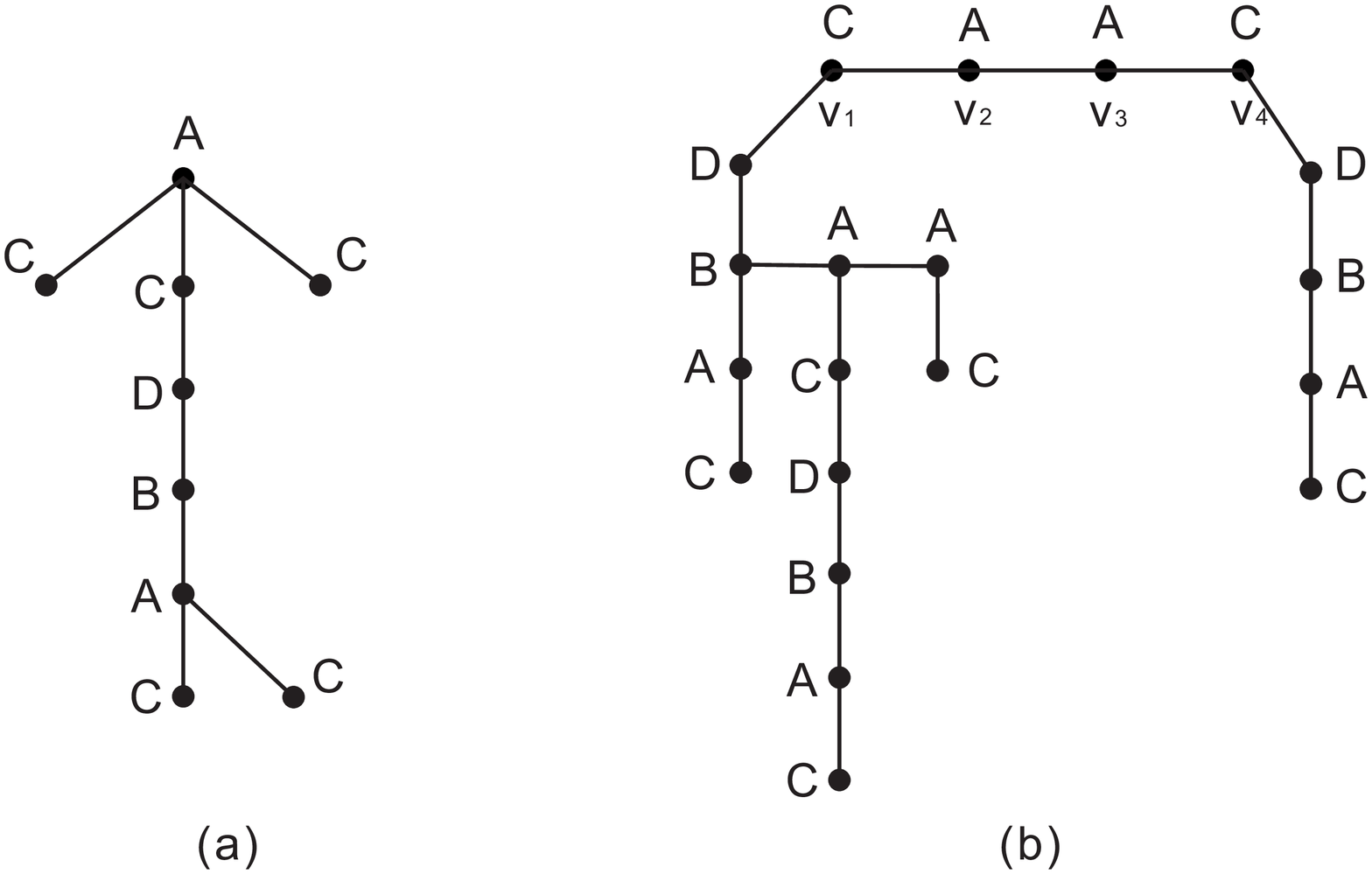}
 \end{center}
\qquad \qquad \qquad \qquad \qquad \qquad \qquad \qquad \qquad
{\small Fig.2}$\\$

Let $(T, S)\in \mathscr{T}_1$ (or $\mathscr{T}_2$) be a labeled tree
for some labeling $S$. Then there is a sequence of labeled trees
$(T_0, S_0)$, $(T_1, S_1), \cdots, (T_{k-1}, S_{k-1})$, $(T_k, S_k)$
such that $(T_0, S_0)=(P_3, S')$ (or $(P_4, S'')$), $(T_k, S_k)=(T,
S)$. The labeled tree $(T_i, S_i)$ can be obtained from $(T_{i-1},
S_{i-1})$ by one of the operations $\mathscr{O}_1$ and
$\mathscr{O}_3$ (or $\mathscr{O}_2$, $\mathscr{O}_3$ and
$\mathscr{O}_4$), where $i\in \{1, 2, \cdots, k\}$. We call the
number of terms in such a sequence of labeled trees that is used to
construct $(T, S)$, the \emph{length} of the sequence. Clearly, the
above sequence has length $k$. We remark that a sequence of labeled
trees used to construct $(T, S)$ is not necessarily unique.


Two main conclusions of our paper are listed as follows.

\begin{thm}
If $T$ is a nontrivial tree of order $n(T)$ with $l(T)$ leaves, then
$\gamma^{d}_2(T)\geq \frac{n(T)-l(T)+3}{4}$, with equality if and
only if $(T, S)\in \mathscr{T}_1$ for some labeling $S$.
\end{thm}

\begin{thm}
If $T$ is a nontrivial tree of order $n(T)$ with $l(T)$ leaves and
$s(T)$ support vertices, then $\gamma^{d}_2(T)\leq
\frac{n(T)+l(T)+s(T)}{4}$, with equality if and only if $(T, S)\in
\mathscr{T}_2$ for some labeling $S$.
\end{thm}

Furthermore, we can slightly improve the upper bound of Theorem~2.4.

\begin{cor}
If $T$ is a nontrivial tree of order $n(T)$ with $l(T)$ leaves and
$s(T)$ support vertices, then $\gamma^{d}_2(T)\leq
\frac{n(T)+3s(T)-l(T)}{4}$.
\end{cor}

\bp Let $T'$ be the tree obtained from $T$ by deleting all but one
leaf from each support vertex of $T$. Then,
$n(T')=n(T)-[l(T)-s(T)]$, $s(T')=s(T)$, $l(T')=s(T)$ and
$\gamma^{d}_2(T)=\gamma^{d}_2(T')$. By Theorem~2.4, we have that
$\gamma^{d}_2(T)=\gamma^{d}_2(T')\leq \frac{n(T')+l(T')+s(T')}{4}=
\frac{n(T)-[l(T)-s(T)]+2s(T)]}{4}=\frac{n(T)+3s(T)-l(T)}{4}$. \ep

\section{Proof of Theorem~2.3}
The following observation establishes properties of trees in the
family $\mathscr{T}_1$.

\begin{obs}
If $(T, S)\in \mathscr{T}_1$, then $(T, S)$ has the following
properties.

$(a)$ Every support vertex of $T$ has status $A$ and every leaf has
status $C$.

$(b)$ Let $v$ be a vertex has status $A$, then sta$(u)\in \{B, C\}$
for $u\in N(v)$.

$(c)$ The set $S_A$ is a $2DD$-set of $T$.

$(d)$ The set $S_A$, $S_B$, $S_C$ and $S_D$ are independent sets.

$(e)$ If sta$(v)\neq A$, then $d(v)\leq 2$.
\end{obs}

\begin{lem}
If $T$ is a tree of order $n(T)\geq 3$ with $l(T)$ leaves, and $(T,
S)\in \mathscr{T}_1$ for some labeling $S$, then
$\gamma^{d}_2(T)=|S_A|=\frac{n(T)-l(T)+3}{4}$, and the set $S_A$ is
the unique $\gamma^{d}_2$-set of $T$.
\end{lem}

\bp We proceed by induction on the length $k$ of a sequence required
to construct the labeled tree $(T, S)$. Let $D$ be any
$\gamma^{d}_2$-set of $T$.

When $k=0$, $(T, S)=(P_3, S')$, $\gamma^{d}_2(T)=|S_A|=1$, the set
$S_A$ is the unique $\gamma^{d}_2$-set of $T$. This establishes the
base case. Let $k\geq 1$ and assume that if the length of sequence
used to construct a labeled tree $(T', S^{*})\in \mathscr{T}_1$ is
less than $k$, then
$\gamma^{d}_2(T')=|S^{*}_A|=\frac{n(T')-l(T')+3}{4}$, $S^{*}_A$ is
the unique $\gamma^{d}_2$-set of $T'$. Now, $(T, S)\in
\mathscr{T}_1$ and let $(T_0, S_0)$, $(T_1, S_1), \cdots, (T_{k-1},
S_{k-1})$, $(T_k, S_k)$ be a sequence of length $k$ used to
construct $(T, S)$, where $(T_0, S_0)=(P_3, S')$, $(T_k, S_k)=(T,
S)$, $(T_i, S_i)$ can be obtained from $(T_{i-1}, S_{i-1})$ by one
of the operations $\mathscr{O}_1$ and $\mathscr{O}_3$, $i\in \{1, 2,
\cdots, k\}$. Let $T'=T_{k-1}$ and $S^{*}=S_{k-1}$. Note that $(T',
S^{*})\in \mathscr{T}_1$. By the inductive hypothesis,
$\gamma^{d}_2(T')=|S^{*}_A|=\frac{n(T')-l(T')+3}{4}$, $S^{*}_A$ is
the unique $\gamma^{d}_2$-set of $T'$. $(T, S)$ can be obtained from
$(T', S^{*})$ by operation $\mathscr{O}_1$ or $\mathscr{O}_3$.

In the former case, we have that $n(T)=n(T')+1$, $l(T)=l(T')+1$, and
$|S_A|=|S^{*}_A|$. It follows Observation~3.1(c) that
$\gamma^{d}_2(T)\leq
|S_A|=|S^{*}_A|=\frac{n(T')-l(T')+3}{4}=\frac{n(T)-1-l(T)+1+3}{4}=\frac{n(T)-l(T)+3}{4}$.
On the other hand, assume that $V(T)\setminus V(T')=\{u\}$, and $v$
is the support vertex of $u$. Take a set $D'=(D\setminus (L(T)\cap
N(v)))\cup \{v\}$ when $(L(T)\cap N(v))\cap D\neq \emptyset$,
otherwise, $D'=D$. $D'$ is a $2DD$-set of $T'$. That is,
$\gamma^{d}_2(T)\geq \gamma^{d}_2(T')=|S^{*}_A|=|S_A|$. In summary,
$\gamma^{d}_2(T)=|S_A|=\frac{n(T)-l(T)+3}{4}$. By the inductive
hypothesis, $S^{*}_A$ is the unique $\gamma^{d}_2$-set of $T'$.
Hence, $D'=S^{*}_A$. In addition, if $u\in D$, then $v\not \in D$.
It follows from $(T, S)\in \mathscr{T}_1$ and Observation~3.1(a),
(b) that $v$ has status $A$, and the non-leaf neighbor of $v$, say
$w$, has status $B$ or $C$. From the choice of $D'$ and
$D'=S^{*}_A$, $u$ is the unique vertex in $D$ which is within
distance two from $w$. It conclude that $w$ is not $2D$-dominated by
$D$, a contradiction. Therefore, $u\not \in D$. Similarly, all
leaf-neighbors of $v$ do not belong to $D$, and then
$D=D'=S^{*}_A=S_A$.

In the latter case, the tree $T$ obtained from $T'$ by attaching a
path $P_4=u_1u_2u_3u_4$ to a leaf $v$ of $T'$, where $u_4$ is a leaf
in $T$. We have that $n(T)=n(T')+4$, $l(T)=l(T')$ and
$|S_A|=|S^{*}_A|+1$. It follows Observation~3.1(c) that
$\gamma^{d}_2(T)\leq
|S_A|=|S^{*}_A|+1=\frac{n(T')-l(T')+3}{4}+1=\frac{n(T)-4-l(T)+3}{4}+1=\frac{n(T)-l(T)+3}{4}$.
Let $D'=(D\setminus \{u_4\})\cup \{u_3\}$ when $u_4\in D$ and $D'=D$
when $u_4\not \in D$, $D''=(D'\setminus \{u_1, u_2\})\cup \{v\}$
when $u_1$ or $u_2$ belong to $D'$, otherwise, $D''=D'$. Then
$u_3\in D$ and $D''\setminus \{u_3\}$ is a $2DD$-set of $T'$. That
is, $\gamma^{d}_2(T)-1\geq \gamma^{d}_2(T')=|S^{*}_A|=|S_A|-1$. In
summary, $\gamma^{d}_2(T)=|S_A|=\frac{n(T)-l(T)+3}{4}$. By the
inductive hypothesis, $S^{*}_A$ is the unique $\gamma^{d}_2$-set of
$T'$. Hence, $D''\setminus \{u_3\}=S^{*}_A$. If $|\{u_1, u_2, u_3,
u_4, v\}\cap D|\geq 2$, the set $(D\setminus \{u_1, u_2, u_3,
u_4\})\cup \{v\}$ is a $2DD$-set of $T'$. More precisely,
$(D\setminus \{u_1, u_2, u_3, u_4\})\cup \{v\}$ is a
$\gamma^{d}_2$-set of $T'$. By the uniqueness of $\gamma^{d}_2$-set
of $T'$, $(D\setminus \{u_1, u_2, u_3, u_4\})\cup \{v\}=S^{*}_A$, a
contradiction. Hence, $|\{u_1, u_2, u_3, u_4, v\}\cap D|=1$. It
implies that $\{u_1, u_2, u_3, u_4, v\}\cap D=\{u_3\}$. It is easy
to see that $D\setminus \{u_3\}$ is a $\gamma^{d}_2$-set of $T'$. By
the uniqueness of $\gamma^{d}_2$-set of $T'$, $D\setminus
\{u_3\}=S^{*}_A$. So, $D=S_A$.
 \ep

In what follows, we begin to prove Theorem~2.3.

\bp The sufficiency follows immediately from Lemma~3.2. So we prove
the necessity only. If $diam(T)\leq 2$, $T$ is a star,
$\gamma^{d}_2(T)=1\geq \frac{n(T)-l(T)+3}{4}$. Suppose that
$\gamma^{d}_2(T)=\frac{n(T)-l(T)+3}{4}$, it is easy to see that
there exists a labeling $S$ of the vertices of $T$ such that $(T,
S)$ can be obtained from $(P_3, S')$ by repeated applications of
operation $\mathscr{O}_1$. Hence, $(T, S)\in \mathscr{T}_1$. If
$diam(T)=3$, $T$ is a double star, and then
$\gamma^{d}_2(T)=2>\frac{n(T)-l(T)+3}{4}$. So, we assume that
$diam(T)\geq 4$. The proof is by induction on $n(T)$. The result is
immediate for $n(T)\leq 5$. For the inductive hypothesis, let
$n(T)\geq 6$. Assume that for every nontrivial tree $T'$ of order
less than $n(T)$, we have that $\gamma^{d}_2(T')\geq
\frac{n(T')-l(T')+3}{4}$, with equality only if $(T', S^{*})\in
\mathscr{T}_1$ for some labeling $S^{*}$.

Let $D$ be a $\gamma^{d}_2$-set of $T$ which contains no leaf and
$P=v_1v_2\cdots v_t$ be a longest path in $T$ such that $d(v_3)$ as
large as possible.

We now proceed with a series of claims that we may assume are
satisfied by the tree $T$, for otherwise the desired result holds.

{\flushleft\textbf{Claim 1.}}\quad Each support vertex in $T$ has
exactly one leaf-neighbor.

If not, assume that there is a support vertex $u$ which is adjacent
to at least two leaves. Deleting one of its leaf-neighbors, say
$u_1$, and denote the resulting tree by $T'$. Observe that
$n(T)=n(T')+1$, $l(T)=l(T')+1$ and $D$ is still a $2DD$-set of $T'$.
That is, $\gamma^{d}_2(T)\geq \gamma^{d}_2(T')\geq
\frac{n(T')-l(T')+3}{4}=\frac{n(T)-1-l(T)+1+3}{4}=\frac{n(T)-l(T)+3}{4}$.

In particular, if $\gamma^{d}_2(T)=\frac{n(T)-l(T)+3}{4}$, then
$\gamma^{d}_2(T')=\frac{n(T')-l(T')+3}{4}$. It means that $(T',
S^{*})\in \mathscr{T}_1$ for some labeling $S^{*}$. By
Observation~3.1(a), $u$ has status $A$. Let $S$ be obtained from
$S^{*}$ by labeling $u_1$ with label $C$. Then $(T, S)$ can be
obtained from $(T', S^{*})$ by operation $\mathscr{O}_1$. Thus, $(T,
S)\in \mathscr{T}_1$. \ep

By Claim~1, we can assume that $d(v_2)=2$. And by Corollary~2.2,
$v_2\in D$. Now, we consider the vertex $v_3$.

{\flushleft\textbf{Claim 2.}}\quad $d(v_3)=2$.

Suppose that $d(v_3)\geq 3$. If $v_3\in D$, let $T'=T-\{v_1, v_2\}$.
Clearly, $D\setminus \{v_2\}$ is a $2DD$-set of $T'$. Note that
$n(T)=n(T')+2$, $l(T)=l(T')+1$, then $\gamma^{d}_2(T)\geq
\gamma^{d}_2(T')+1\geq
\frac{n(T')-l(T')+3}{4}+1=\frac{n(T)-2-l(T)+1+3}{4}+1>\frac{n(T)-l(T)+3}{4}$.
So we assume that $v_3\not \in D$. If $v_3$ is adjacent to a support
vertex outside $P$, say $v_2'$. It follows from Claim~1 and
Corollary~2.2 that $v_2'\in D$. Moreover, $(D\setminus \{v_2,
v_2'\})\cup \{v_3\}$ is a $2DD$-set of the tree $T'$ obtained from
$T$ by removing all leaf-neighbors of $v_2$ and $v_2'$. Hence,
$\gamma^{d}_2(T)\geq \gamma^{d}_2(T')+1\geq
\frac{n(T')-l(T')+3}{4}+1=\frac{n(T)-2-l(T)+3}{4}+1>\frac{n(T)-l(T)+3}{4}$.
Combining the assumption that $d(v_3)\geq 3$, $v_3$ is a support
vertex of degree three of $T$. We remove its leaf-neighbor, say $u$,
and $D$ is still a $2DD$-set of the resulting tree $T'$ from $u\not
\in D$. Hence, $\gamma^{d}_2(T)\geq \gamma^{d}_2(T')\geq
\frac{n(T')-l(T')+3}{4}=\frac{n(T)-l(T)+3}{4}$. We show that in fact
$\gamma^{d}_2(T)>\frac{n(T)-l(T)+3}{4}$. Suppose to the contrary
that $\gamma^{d}_2(T)=\frac{n(T)-l(T)+3}{4}$. Then we have equality
throughout the above inequality chain. In particular,
$\gamma^{d}_2(T)=\gamma^{d}_2(T')=\frac{n(T')-l(T')+3}{4}$. By the
inductive hypothesis, $(T', S^{*})\in \mathscr{T}_1$ for some
labeling $S^{*}$. By Observation~3.1(a) and (b), the vertex $v_3$
has status $B$ or $C$ in $S^{*}$. Since $D$ contains no leaf, $D$ is
also a $\gamma^{d}_2$-set of $T'$. On the other hand, by Lemma~3.2,
$S^{*}_A$ is the unique $\gamma^{d}_2$-set of $T'$. So, $D=S^{*}_A$.
It implies that $u$ can not be $2D$-dominated by $D$, a
contradiction.
 \ep

{\flushleft\textbf{Claim 3.}}\quad $d(v_4)=2$.

Assume that $d(v_4)\geq 3$ and $v_3'$ is a neighbor of $v_4$ outside
$P$. From Claim~1 and the choice of $P$, one of the three cases as
following holds:

(1) $v_3'$ is adjacent to a support vertex, say $v_2'$, where $v_2'$
and $v_3'$ have degree two;

(2) $v_3'$ is a support vertex of degree two in $T$;

(3) $v_3'$ is a leaf.

In the first case, let $T'$ be a tree obtained from $T$ by removing
$v_1, v_2, v_3$ and the leaf-neighbor of $v_2'$. We have that
$n(T)=n(T')+4$, $l(T)=l(T')+1$ and $\gamma^{d}_2(T')\leq
\gamma^{d}_2(T)-1$. In the latter two cases, let $T'=T-\{v_1, v_2,
v_3\}$. We have that $n(T)=n(T')+3$, $l(T)=l(T')+1$ and
$\gamma^{d}_2(T')\leq \gamma^{d}_2(T)-1$. In either case, we always
have $\gamma^{d}_2(T)>\frac{n(T)-l(T)+3}{4}$ by an argument similar
to the proof of Claim~2. \ep

Let $T'=T-\{v_1, v_2, v_3, v_4\}$. Note that $n(T)=n(T')+4$,
$\gamma^{d}_2(T')\leq \gamma^{d}_2(T)-1$. In addition,
$l(T)=l(T')+1$ when $d(v_5)\geq 3$, and $l(T)=l(T')$ when
$d(v_5)=2$. Hence, we always have that $\gamma^{d}_2(T)\geq
\gamma^{d}_2(T')+1\geq \frac{n(T')-l(T')+3}{4}+1\geq
\frac{n(T)-4-l(T)+3}{4}+1=\frac{n(T)-l(T)+3}{4}$. Suppose that
$\gamma^{d}_2(T)=\frac{n(T)-l(T)+3}{4}$, then we have equality
throughout the above inequality chain. In particular, $d(v_5)=2$ and
$\gamma^{d}_2(T)-1=\gamma^{d}_2(T')=\frac{n(T')-l(T')+3}{4}$. By the
inductive hypothesis, $(T', S^{*})\in \mathscr{T}_1$ for some
labeling $S^{*}$. Since $v_5$ is a leaf in $T'$, by
Observation~3.1(a), it has status $C$. Let $S$ be obtained from the
labeling $S^{*}$ by labeling the vertices $v_1, v_2, v_3, v_4$ with
label $C, A, B, D$, respectively. Then, $(T, S)$ can be obtained
from $(T', S^{*})$ by operation $\mathscr{O}_3$. Thus, $(T, S)\in
\mathscr{T}_1$. \ep

\section{Proof of Theorem~2.4}

The following observation establishes properties of trees in the
family $\mathscr{T}_2$.

\begin{obs}
If $(T, S)\in \mathscr{T}_2$, then $(T, S)$ has the following
properties.

$(a)$ Every support vertex of $T$ has status $A$ and every leaf has
status $C$.

$(b)$ The set $S_A$ is a $2DD$-set of $T$.

$(c)$ Let $v$ be a vertex which has status $A$ or $B$, $v$ has at
most one corresponding vertex. In particular, if there is no
corresponding vertex of degree two of $v$ in $T$, then $d(v)=2$.

$(d)$ If $v$ is a support vertex, then $v$ has degree two.

$(e)$ Let $v$ be a vertex of degree two which has status $C$, then
it is adjacent to two vertices, say $u$ and $w$, which are labeled
$A$ and $D$, respectively. In particular, if $d(u)=2$, the component
of $T-vw$ containing $v$, say $T'$, containing the basic path of
$T$, and $(T', S^{*})\in \mathscr{T}_2$ for some labeling $S^{*}$.
\end{obs}

\begin{lem}
Let $T$ be a tree and $S$ be a labeling of $T$ such that $(T, S)\in
\mathscr{T}_2$. Then, $\gamma^{d}_2(T)=\frac{n(T)+s(T)+l(T)}{4}$.
\end{lem}

\bp By Observation~4.1(b), $S_A$ is a $2DD$-set of $T$ and
$S_A=\frac{n(T)+s(T)+l(T)}{4}$ (We can obtain this conclusion by
induction on $n(T)$, it is similar to the proof of Lemma~3.2, so we
omit it). So, $\gamma^{d}_2(T)\leq \frac{n(T)+s(T)+l(T)}{4}$. Since
$(T, S)\in \mathscr{T}_2$, $T=P_4$ when $n\leq 4$, and
$\gamma^{d}_2(T)=2=\frac{n(T)+s(T)+l(T)}{4}$. So, we assume that
$n(T)\geq 5$. Combining the definition of $\mathscr{T}_2$, we have
that $diam(T)\geq 7$. Suppose that $T$ is a tree with minimum order
which satisfy the two properties:

(1) $(T, S)\in \mathscr{T}_2$;

(2) $\gamma^{d}_2(T)<\frac{n(T)+s(T)+l(T)}{4}$.

Let $D$ be a $\gamma^{d}_2$-set of $T$ which contains no leaf,
$u_1u_2u_3u_4$ be the basic path of $T$, and $v_1$ be a leaf of $T$
that at maximum distance from $u_2$, let $P=v_1v_2v_3\cdots v_tu_2$
be the path between $v_1$ and $u_2$. Note that $v_t=u_1$ or $u_3$.
It follows from $(T, S)\in \mathscr{T}_2$ and Observation~4.1(d)
that $d(v_2)=2$ and $v_1, v_2$ have status $C, A$, respectively. And
moreover, by the definition of $\mathscr{T}_2$, $v_3$ has status $A$
or $B$.

In the form case, if $d(v_3)=2$, then $v_1v_2v_3v_4$ is the basic
path of $T$, a contradiction. So, $d(v_3)\geq 3$. It implies that
there exists a sequence of length $k$ used to construct $(T, S)$:
$(P_4, S'')$, $(T_1, S_1), \cdots, (T_{k-1}, S_{k-1})$, $(T, S)$,
such that $(T, S)$ is obtained from $(T_{k-1}, S_{k-1})$ by
operation $\mathscr{O}_4$. That is, $T$ is obtained from $T_{k-1}$
by adding the path $v_1v_2$ and joining $v_2$ to $v_3$. But in this
case, by the definition of $\mathscr{O}_4$, we can always obtain a
leaf which is farther away from $u_2$ than $v_1$, contradicting the
choice of $v_1$. So we assume that $v_3$ has status $B$.

If $d(v_3)\geq 3$, by Observation~4.1(d), $v_3$ is not a support
vertex. From the choice of $v_1$ and the fact that $diam(T)\geq 7$,
$v_3$ is adjacent to $s$ support vertices of degree two other than
$v_2$, where $s\geq 1$. These support vertices are labeled $A$, and
the leaf-neighbor of each of them is labeled $C$. From the choice of
$D$ and Corollary~2.2, $S(T)\cap N(v_3)\subseteq D$.
 $v_4, v_5, v_6$ has status $D, C, A$, respectively, and $d(v_4)=d(v_5)=2$.
  Moreover, there exists no a corresponding vertex of degree two of $v_6$ in $T$, so $d(v_6)=2$.
 Note that $\{v_3,
v_4, v_5, v_6\}\cap D\neq \emptyset$, then $(D\setminus \{v_3, v_4,
v_5\})\cup \{v_6\}$ is also a $\gamma^{d}_2$-set of $T$. Hence,
$D'=D\setminus \{v_2\}$ is a $2DD$-set of $T'$ with order at most
$\gamma^{d}_2(T)-1$, where $T'=T-\{v_1, v_2\}$. On the other hand,
note that $(T', S^{*})\in \mathscr{T}_2$ for some labeling $S^{*}$,
from the choice of $T$,
$\gamma^{d}_2(T')=\frac{n(T')+s(T')+l(T')}{4}=\frac{n(T)+s(T)+l(T)}{4}-1>\gamma^{d}_2(T)-1$.
A contradiction.

If $d(v_3)=2$, from the definition of $\mathscr{T}_2$, $v_4$ has
status $D$, and furthermore, $v_5, v_6$ have status $C, A$,
respectively. In particular, $d(v_4)=d(v_5)=2$. Note that $v_2\in
D$, and $\{v_3, v_4, v_5, v_6\}\cap D\neq \emptyset$, so the set
$D'=(D\setminus \{v_3, v_4, v_5\})\cup \{v_6\}$ is also a
$\gamma^{d}_2$-set of $T$. Now, we distinguish two cases as follows.

{\flushleft\textbf{Case 1.}}\quad $d(v_6)=2$.

The set $D''=D'\setminus \{v_2\}$ is a $2DD$-set of $T'$ with order
at most $\gamma^{d}_2(T)-1$, where $T'=T-\{v_1, v_2, v_3, v_4\}$. On
the other hand, from the choice of $T$ and the fact that $(T',
S^{*})\in \mathscr{T}_2$ for some labeling $S^{*}$,
$\gamma^{d}_2(T')=\frac{n(T')+s(T')+l(T')}{4}=\frac{n(T)+s(T)+l(T)}{4}-1>\gamma^{d}_2(T)-1$.
A contradiction.

{\flushleft\textbf{Case 2.}}\quad $d(v_6)\geq 3$.

We have that sta$(v_7)=A$ or $B$. If sta$(v_7)=B$, then all
neighbors of $v_6$ outside $P$ have status $A$, and note that these
neighbors are support vertices of degree two (From the choice of
$v_1$ and the definition of $\mathscr{T}_2$). We remove one of these
support vertices, say $u_1$, and its leaf-neighbor, say $u_2$,
denote the resulting tree by $T'$. Clearly, $(T', S^{*})\in
\mathscr{T}_2$ for some labeling $S^{*}$. We know that $v_2, v_6\in
D'$, and $\{u_1, u_2\}\cap D'\neq \emptyset$, so $D''=D'\setminus
\{u_1, u_2\}$ is a $2DD$-set of $T'$ with order at most
$\gamma^{d}_2(T)-1$, where $T'=T-\{u_1, u_2\}$. On the other hand,
from the choice of $T$,
$\gamma^{d}_2(T')=\frac{n(T')+s(T')+l(T')}{4}=\frac{n(T)+s(T)+l(T)}{4}-1>\gamma^{d}_2(T)-1$.
A contradiction.

If sta$(v_7)=A$, then one of the two cases as following holds:

(1) There exists a neighbor of $v_6$ outside $P$, say $u_1$, has
status $B$.

(2) All neighbors of $v_6$ outside $P$ have status $A$.

In the former case, there exists a neighbor $u_2$ of $u_1$ which has
status $D$. Similarly, there exists a neighbor $u_3$ of $u_2$ which
has status $C$, and there exists a neighbor $u_4$ of $u_3$ which has
status $A$. Moreover, let $u_5$ be a neighbor of $u_4$ other than
$u_3$, then $u_5$ has status $A$ or $B$. In either case, $u_5$ has
degree at least two, which contradicts the choice of $v_1$.

In the latter case, we take any neighbor of $v_6$ outside $P$, say
$u_1$, and we have that $u_1$ has a neighbor which has status $C$,
say $u_2$. From the choice of $v_1$, $u_2$ is a leaf. By
Observation~4.1(d), $d(u_1)=2$. And we can obtain a contradiction by
an argument similar to the case that sta$(v_7)=B$ as above.

In summary, if $(T, S)\in \mathscr{T}_2$. Then,
$\gamma^{d}_2(T)=\frac{n(T)+s(T)+l(T)}{4}$. \ep

\begin{lem}
Let $T$ be a tree and $S$ be a labeling of $T$ such that $(T, S)\in
\mathscr{T}_2$. Then for any leaf $v$, there exists a set $D$ with
order $\frac{n(T)+s(T)+l(T)}{4}-1$ such that each vertex of $T$ is
$2D$-dominated by $D$ except for $v$, and the non-leaf neighbor of
the support vertex of $v$ belongs to $D$.
\end{lem}

\bp Take any leaf $v_1$ of $T$. We proceed by induction on the
length $k$ of a sequence required to construct the labeled tree $(T,
S)$. When $k=0$, $(T, S)=(P_4, S'')$, the result is immediate. Let
$k\geq 1$ and assume that if the length of sequence used to
construct a labeled tree $(T', S^{*})\in \mathscr{T}_2$ is less than
$k$, the result holds. Since $(T, S)\in \mathscr{T}_2$, there exists
always a sequence of length $k$ used to construct $(T, S)$: $(P_4,
S'')$, $(T_1, S_1), \cdots, (T_{k-1}, S_{k-1})$, $(T, S)$.

First, we assume that $v_1$ is in the basic path of $T$. Since
$(T_{k-1}, S_{k-1})\in \mathscr{T}_2$, $v_1$ is still a leaf of
$T_{k-1}$. By the inductive hypothesis, there exists a set $D'$ with
order $\frac{n(T_{k-1})+s(T_{k-1})+l(T_{k-1})}{4}-1$ such that each
vertex of $T_{k-1}$ is $2D$-dominated by $D'$ except for $v_1$, and
$v_3$ belongs to $D'$, where $v_3$ is the neighbor of the support
vertex of $v_1$. We know that $(T, S)$ is obtained from $(T_{k-1},
S_{k-1})$ by one of the operations $\mathscr{O}_2$, $\mathscr{O}_3$
and $\mathscr{O}_4$. In the first or third case, let $D$ be the set
consisting of $D'$ and the support vertex which belongs to
$V(T)\setminus V(T_{k-1})$, and $D$ is the desired set. In the
second case, the tree $T$ is obtained from $T_{k-1}$ by adding a
path $u_1u_2u_3u_4$ and joining $u_1$ to a leaf $u$ of $T_{k-1}$.
Note that $u$ has status $C$, and by Observation~4.1(d), the
neighbor of $u$, say $u'$, has degree two. By the inductive
hypothesis, there exists a set $D'$ with order
$\frac{n(T_{k-1})+s(T_{k-1})+l(T_{k-1})}{4}-1$ such that each vertex
of $T_{k-1}$ is $2D$-dominated by $D'$ except for $v_1$, and $v_3$
belongs to $D'$. Moreover, one of $u$ and $u'$ belongs to $D'$. Let
$D$ be the set consisting of $D'$ and the vertex $u_3$, and $D$ is
the desired set.

Next, we consider the case that $v_1$ is not in the basic path.
Since $(T, S)\in \mathscr{T}_2$, this leaf has status $C$ and its
support vertex $v_2$ is labeled $A$. By Observation~4.1(d), $v_2$
has degree two. Let $P=v_1v_2\cdots v_tv$ be the path between $v_1$
and $v$, where $v$ is the vertex of basic path which has minimum
distance from $v_1$. Note that the neighbor of $v_2$, say $v_3$, has
status $A$ or $B$.


Next, we distinguish two cases as follows.

{\flushleft\textbf{Case 1.}}\quad sta$(v_3)=A$.

If $d(v_3)=2$, then it is easy to see that $v_1v_2v_3v_4$ is the
basic path of $T$, a contradiction.

If $d(v_3)\geq 3$, from the definition of $\mathscr{T}_2$ and the
fact that a sequence of labeled trees used to construct $(T, S)$ is
not necessarily unique, we have that there exists a sequence of
length $k$ used to construct $(T, S)$: $(P_4, S'')$, $(T_1', S_1'),
\cdots, (T_{k-1}', S_{k-1}')$, $(T, S)$, such that $(T, S)$ is
obtained from $(T_{k-1}', S_{k-1}')$ by $\mathscr{O}_4$. That is,
the tree $T$ is obtained from $T_{k-1}'$ by adding the path $v_1v_2$
and joining $v_2$ to a vertex $v_3$. Note that $v_3$ has a neighbor
of degree two, say $u$, which is labeled $C$ (Otherwise, no vertex
of $T$ is the corresponding vertex of $v_3$). By Observation~4.1(e),
the component of $T_{k-1}'-uu'$ containing $u$, say $T'$, containing
the basic path, and $(T', S^{*})\in \mathscr{T}_2$ for some $S^{*}$,
where $u'$ is the neighbor of $u$ other than $v_3$. It implies that
there always exists a sequence of length $k$ used to construct $(T,
S)$: $(P_4, S'')$, $(T_1'', S_1''), \cdots, (T_{k-1}'', S_{k-1}'')$,
$(T, S)$, satisfying the two conditions as follows:

(1) $(T_{k-1}'', S_{k-1}'')=(T_{k-1}', S_{k-1}')$;

(2) There is a $i\in \{1, 2, \cdots, k-2\}$ in this sequence such
that $(T_i'', S_i'')=(T', S^{*})$.

By the inductive hypothesis, there exists a set $D'$ with order
$\frac{n(T_i'')+s(T_i'')+l(T_i'')}{4}-1$ such that each vertex of
$T_i''$ is $2D$-dominated by $D'$ except for $u$, and $u''$ belongs
to $D'$, where $u''$ is a neighbor of $v_3$ in $T_i''$ other than
$u$. Then $D_i=D'\cup \{v_3\}$ is a $\gamma^{d}_2$-set of $T_i''$.
For each $j\in \{i, i+1, \cdots, k-2\}$, we know that $(T_{j+1}'',
S_{j+1}'')$ is obtained from $(T_j'', S_j'')$ by one of the
operations $\mathscr{O}_2$, $\mathscr{O}_3$ and $\mathscr{O}_4$. Let
$D_{j+1}=D_j\cup \{w\}$, where $w\in V(T_{j+1}'')\setminus V(T_j)''$
and has status $A$. It is easy to see that $D_{j+1}$ is a
$\gamma^{d}_2$-set of $T_{j+1}''$, and moreover, $D_{k-1}$ is the
desired set.



{\flushleft\textbf{Case 2.}}\quad sta$(v_3)=B$.\

In this case, if $d(v_3)\geq 3$, there must be a neighbor of $v_3$,
say $u$, which has status $A$. From the definition of
$\mathscr{T}_2$, the component of $T-v_3u$ containing $v_3$, say
$T'$, containing the basic path, and $(T', S^{*})\in \mathscr{T}_2$
for some $S^{*}$. We can obtain the desired set by an argument
similar to the case of sta$(v_3)=A$ and $d(v_3)\geq 3$.

If $d(v_3)=2$, then $v_4, v_5, v_6$ have status $D, C, A$,
respectively, and $d(v_4)=d(v_5)=2$. If $d(v_6)=2$, let $T'=T-\{v_1,
v_2, v_3, v_4\}$. Note that $(T', S^{*})\in \mathscr{T}_2$ for some
$S^{*}$. By the inductive hypothesis, there exists a set $D'$ with
order $\frac{n(T')+s(T')+l(T')}{4}-1$ such that each vertex of $T'$
is $2D$-dominated by $D'$ except for $v_5$, and $v_7$ belongs to
$D'$, then the set $D'\cup \{v_3\}$ is the desired set. So we
consider the case of $d(v_6)\geq 3$. From the definition of
$\mathscr{T}_2$, there must exist a neighbor of $v_6$, say $u$, such
that sta$(u)=A$ and the component of $T-v_6u$ containing $v_6$, say
$T'$, containing the basic path, and $(T', S^{*})\in \mathscr{T}_2$
for some $S^{*}$. By the inductive hypothesis, there exists a set
$D'$ with order $\frac{n(T')+s(T')+l(T')}{4}-1$ such that each
vertex of $T'$ is $2D$-dominated by $D'$ except for $v_1$, and $v_3$
belongs to $D'$, We can obtain the desired set by an argument
similar to the case of sta$(v_3)=A$ and $d(v_3)\geq 3$. \ep

In what follows, we begin to prove Theorem~2.3.

\bp The sufficiency follows immediately from Lemma~4.2. So we prove
the necessity only. If $diam(T)\leq 2$, $T$ is a star, and
$\gamma^{d}_2(T)=1<\frac{n(T)+s(T)+l(T)}{4}$. If $diam(T)=3$, $T$ is
a double star, and then $\gamma^{d}_2(T)=2\leq
\frac{n(T)+s(T)+l(T)}{4}$. Support that $\gamma^{d}_2(T)=
\frac{n(T)+s(T)+l(T)}{4}$, it is easy to see that $T=P_4$, let $S$
be the labeling that assigns to the two leaves of the path $P_4$
status $C$, and the remaining vertices status $A$, then the label
tree $(P_4, S)\in \mathscr{T}_2$. So we assume that $diam(T)\geq 4$.
The proof is by induction on $n(T)$. The result is immediate for
$n(T)\leq 4$. For the inductive hypothesis, let $n(T)\geq 5$. Assume
that for every nontrivial tree $T'$ of order less than $n(T)$, we
have that $\gamma^{d}_2(T')\leq \frac{n(T')+s(T')+l(T')}{4}$, with
equality only if $(T', S^{*})\in \mathscr{T}_2$ for some labeling
$S^{*}$.

Let $D$ be a $\gamma^{d}_2$-set of $T$ which contains no leaf and
$P=v_1v_2\cdots v_t$ be a longest path in $T$ such that

(i) $d(v_5)$ as large as possible, and subject to this condition

(ii) $d(v_4)$ as large as possible, and subject to this condition

(iii) $d(v_3)$ as large as possible.

We now proceed with a series of claims that we may assume are
satisfied by the tree $T$, for otherwise the desired result holds.

{\flushleft\textbf{Claim 1.}}\quad Each support vertex in $T$ has
exactly one leaf-neighbor.

If not, assume that there is a support vertex $u$ which is adjacent
to at least two leaves, say $u_1, u_2$. Deleting $u_1$, and denote
the resulting tree by $T'$. Take a $\gamma^{d}_2$-set of $T'$
contains no leaf, say $D'$. It follows that $u$ is either contained
in $D'$ or has at least two non-leaf neighbors in $D'$, and then
$D'$ is also a $2DD$-set of $T$.  That is, $\gamma^{d}_2(T)\leq
\gamma^{d}_2(T')$. Observe that $n(T)=n(T')+1$, $l(T)=l(T')+1$ and
$s(T)=s(T')$. We have that $\gamma^{d}_2(T)\leq \gamma^{d}_2(T')\leq
\frac{n(T')+s(T')+l(T')}{4}=\frac{n(T)-1+s(T)+l(T)-1}{4}<\frac{n(T)+s(T)+l(T)}{4}$.
 \ep

By Claim~1, we can assume that $d(v_2)=2$. And by Corollary~2.2,
$v_2\in D$. Now, we consider the vertex $v_3$.

{\flushleft\textbf{Claim 2.}}\quad $v_3$ is not a support vertex.

In other words, all neighbors of $v_3$ are support vertices of
degree two, except possibly the vertex $v_4$. If not, support that
$v_3$ is a support vertex and $u$ is the leaf-neighbor. Let
$T'=T-\{v_1, v_2\}$. Note that $n(T)=n(T')+2$, $l(T)=l(T')+1$ and
$s(T)=s(T')+1$, then $\gamma^{d}_2(T)\leq \gamma^{d}_2(T')+1\leq
\frac{n(T')+s(T')+l(T')}{4}+1=\frac{n(T)-2+s(T)-1+l(T)-1}{4}+1=\frac{n(T)+s(T)+l(T)}{4}$.
In particular, if $\gamma^{d}_2(T)=\frac{n(T)+s(T)+l(T)}{4}$, then
$\gamma^{d}_2(T')=\frac{n(T')+s(T')+l(T')}{4}$. It means that $(T',
S^{*})\in \mathscr{T}_2$ for some labeling $S^{*}$. By Lemma~4.3,
there exists a $2DD$-set $S$ of $T'-\{u\}$ with cardinality
$\gamma^{d}_2(T')-1$, and the non-leaf neighbor of $v_3$ in $T'$
belongs to $S$. It is easy to see that $S\cup \{v_2\}$ is a
$2DD$-set of $T$ with cardinality $\gamma^{d}_2(T')$. That is,
$\gamma^{d}_2(T)\leq \gamma^{d}_2(T')$, Contradicting the fact that
$\gamma^{d}_2(T)=\gamma^{d}_2(T')+1$. Hence, we have that
$\gamma^{d}_2(T)<\frac{n(T)+s(T)+l(T)}{4}$. \ep

Let $(S(T)\cap N(v_3))\setminus \{v_4\}=\{w_1, w_2, \cdots, w_t\}$,
where $w_1=v_2$, $t\geq 1$.

{\flushleft\textbf{Claim 3.}}\quad $d(v_4)=2$.

Assume that $d(v_4)\geq 3$, let $T'$ be the component of $T-v_3v_4$
containing $v_4$. It follows from $n(T)=n(T')+1+2t$, $l(T)=l(T')+t$
and $s(T)=s(T')+t$ that $\gamma^{d}_2(T)\leq \gamma^{d}_2(T')+t\leq
\frac{n(T')+s(T')+l(T')}{4}+t=\frac{n(T)-1-2t+s(T)-t+l(T)-t}{4}+t<\frac{n(T)+s(T)+l(T)}{4}$.
\ep

{\flushleft\textbf{Claim 4.}}\quad $d(v_5)=2$.

Assume that $d(v_5)\geq 3$ and $v_4'$ be a neighbor of $v_5$ outside
$P$. If $t=2$, from the choice of $P$ and Claim~1, we only need to
consider the two case as follows (In other cases, let $T'=T-\{v_1,
v_2, v_3, v_4\}$. We can always obtain a $\gamma^{d}_2$-set of $T'$
which contains a vertex $u\in N[v_5]\cap V(T')$. It means that
$\gamma^{d}_2(T)\leq \gamma^{d}_2(T')+1$. Observe that
$n(T)=n(T')+4$, $l(T)=l(T')+1$ and $s(T)=s(T')+1$. We always have
that $\gamma^{d}_2(T)<\frac{n(T)+s(T)+l(T)}{4}$):

(1) $v_5$ is not a support vertex, $v_4'$ is adjacent to a support
vertex $v_3'$, where $v_3'$ and $v_4'$ have degree two.

(2) $v_5$ is not a support vertex and $v_4'$ is adjacent to $h$
support vertices of degree two, where $h\geq 2$.

Let $T'$ is the component of $T-v_5v_4'$ containing $v_5$. In the
former case, $n(T)=n(T')+3$, $l(T)=l(T')+1$, $s(T)=s(T')+1$ and
$\gamma^{d}_2(T)\leq \gamma^{d}_2(T')+1$. In the latter case, note
that it is possible that $v_4'$ is a support vertex, then
$n(T')+2h+1\leq n(T)\leq n(T')+2h+2$, $l(T')+h\leq l(T)\leq
l(T')+h+1$, $s(T')+h\leq s(T)\leq s(T')+h+1$ and
$\gamma^{d}_2(T)\leq \gamma^{d}_2(T')+h$. In either case, we
conclude that $\gamma^{d}_2(T)<\frac{n(T)+s(T)+l(T)}{4}$.

If $t\geq 3$, let $T'$ be the component of $T-v_4v_5$ containing
$v_5$. Observe that $n(T)=n(T')+2+2t$, $l(T)=l(T')+t$ and
$s(T)=s(T')+t$ and $\gamma^{d}_2(T)\leq \gamma^{d}_2(T')+t$.
Analogous to the proof of Case~3, we have that
$\gamma^{d}_2(T)<\frac{n(T)+s(T)+l(T)}{4}$.\ep

{\flushleft\textbf{Claim 5.}}\quad $d(v_6)=2$ or all neighbors of
$v_6$ outside $P$ are support vertices of degree two.

First, we show that $v_6$ is not a support vertex. If not, it
follows from Claim~1 that $v_6$ has one leaf-neighbor, and construct
a tree $T'$ which is obtained from $T$ by removing the leaf-neighbor
of $v_6$ and joining a new vertex to $v_2$. Let $D'$ be a
$\gamma^{d}_2$-set of $T'$ which contains no leaf, then $N(v_3)\cap
S(T)\subseteq D'$. We take a set $D''=(D'\setminus \{v_3, v_4,
v_5\})\cup \{v_6\}$ when $D'\cap \{v_3, v_4, v_5\}\neq \emptyset$,
and otherwise, $D''=D'$. Note that $D''$ is also a $2DD$-set of $T$,
and moreover, $n(T)=n(T')$, $l(T)=l(T')$, $s(T)=s(T')+1$. Hence,
$\gamma^{d}_2(T)\leq \gamma^{d}_2(T')\leq
\frac{n(T')+s(T')+l(T')}{4}=\frac{n(T)+s(T)-1+l(T)}{4}<\frac{n(T)+s(T)+l(T)}{4}$.

Let $u_1$ be a leaf outside $P$ that at maximum distance from $v_6$,
and $P_1=u_1u_2\cdots u_{s-1}u_s$ be the path between $u_1$ and
$v_6$, where $u_s=v_6$. Clearly, $s\leq 6$.

If $s=4$, then we have that $u_3$ is adjacent to $a$ support
vertices of degree two, where $a\geq 1$. Suppose that $u_3$ is not a
support vertex, let $T'$ be the component of $T-u_3v_6$ containing
$v_6$. It follows from $n(T)=n(T')+2a+1$, $l(T)=l(T')+a$ and
$s(T)=s(T')+a$ that $\gamma^{d}_2(T)\leq \gamma^{d}_2(T')+a\leq
\frac{n(T')+s(T')+l(T')}{4}+a=\frac{n(T)-2a-1+s(T)-a+l(T)-a}{4}+a<\frac{n(T)+s(T)+l(T)}{4}$.
So, we assume that $u_3$ has a leaf-neighbor, say $u$, and in this
case, let $T'=T-\{u_1, u_2\}$. Note that $n(T)=n(T')+2$,
$l(T)=l(T')+1$ and $s(T)=s(T')+1$, then $\gamma^{d}_2(T)\leq
\gamma^{d}_2(T')+1\leq
\frac{n(T')+s(T')+l(T')}{4}+1=\frac{n(T)-2+s(T)-1+l(T)-1}{4}+1\leq
\frac{n(T)+s(T)+l(T)}{4}$. In particular, if
$\gamma^{d}_2(T)=\frac{n(T)+s(T)+l(T)}{4}$, then
$\gamma^{d}_2(T')=\frac{n(T')+s(T')+l(T')}{4}$. It means that $(T',
S^{*})\in \mathscr{T}_2$ for some labeling $S^{*}$. By Lemma~4.3,
there exists a $2DD$-set $S$ of $T'-\{u\}$ with cardinality
$\gamma^{d}_2(T')-1$, and a non-leaf neighbor of $u_3$ in $T'$
belongs to $S$. It is easy to see that $S\cup \{u_2\}$ is a
$2DD$-set of $T$ with cardinality $\gamma^{d}_2(T')$. That is,
$\gamma^{d}_2(T)\leq \gamma^{d}_2(T')$, Contradicting the fact that
$\gamma^{d}_2(T)=\gamma^{d}_2(T')+1$.

If $s=5$, by an argument similar to that of Claim~1, Claim~2 and
Claim~3, we have that $d(u_2)=d(u_4)=2$, $u_3$ is not a support
vertex and adjacent to $a$ support vertices of degree two, where
$a\geq 1$. Let $T'$ be the component of $T-u_4v_6$ containing $v_6$
and $D'$ be a $\gamma^{d}_2$-set of $T'$ contains no leaf. If $a\geq
2$, Observe that $D'\cup (S(T)\cap N(u_3))$ is a $2DD$-set of $T$.
Combining the fact that $n(T)=n(T')+2a+2$, $l(T)=l(T')+a$,
$s(T)=s(T')+a$. We have that $\gamma^{d}_2(T)\leq
\gamma^{d}_2(T')+a\leq
\frac{n(T')+s(T')+l(T')}{4}+a=\frac{n(T)-2a-2+s(T)-a+l(T)-a}{4}+a<\frac{n(T)+s(T)+l(T)}{4}$.

So we consider the case of $a=1$. If there is a vertex belonging to
$N[v_6]\cap D'$, then $D'\cup \{u_2\}$ is a $2DD$-set of $T$, and so
$\gamma^{d}_2(T)\leq \gamma^{d}_2(T')+1\leq
\frac{n(T')+s(T')+l(T')}{4}+1=\frac{n(T)-4+s(T)-1+l(T)-1}{4}+1<\frac{n(T)+s(T)+l(T)}{4}$.
So we can assume that $N[v_6]\cap D'=\emptyset$. If $\{v_3,
v_4\}\cap D'\neq \emptyset$, then $(D'\setminus \{v_3, v_4\})\cup
\{v_5\}$ is also a $\gamma^{d}_2$-set of $T'$, and we are done. If
$v_3, v_4\not \in D'$, it follows from $d(v_5)=2$ and $N[v_6]\cap
D'=\emptyset$ that $v_5$ is not $2D$-dominated by $D'$, a
contradiction.

If $s=6$, from Claim~1, Claim~2 and the choice of $T$, we have that
$d(u_2)=d(u_4)=d(u_5)=2$, and $u_3$ is not a support vertex and
adjacent to $a$ support vertices of degree two, where $a\leq t$. Let
$T'$ be the component of $T-v_5v_6$ containing $v_6$ and $D_1$ be a
$\gamma^{d}_2$-set of $T'$ contains no leaf. Note that $S(T)\cap
N(u_3)\subseteq D_1$. Take a set $D'=(D_1\setminus \{u_3, u_4,
u_5\})\cup \{v_6\}$ when $\{u_3, u_4, u_5\}\cap D_1\neq \emptyset$,
and otherwise, $D'=D_1$. Observe that $D'\cup \{w_1, w_2, \cdots,
w_t\}$ is a $2DD$-set of $T$. Combining the fact that
$n(T)=n(T')+2t+3$, $l(T)=l(T')+t$, $s(T)=s(T')+t$. We have that
$\gamma^{d}_2(T)\leq \gamma^{d}_2(T')+t\leq
\frac{n(T')+s(T')+l(T')}{4}+t=\frac{n(T)-2t-3+s(T)-t+l(T)-t}{4}+t<\frac{n(T)+s(T)+l(T)}{4}$.
\ep

We assume that $|N(v_6)\setminus \{v_5, v_7\}|=a$, then $a\geq 0$.
In addition, by the claims as above, we have that
$d(v_2)=d(v_4)=d(v_5)=2$, $v_3$ is not a support vertex and adjacent
to $t$ support vertices of degree two, where $t\geq 1$.

If $a=0$, then $d(v_6)=2$. Let $T'$ be the component of $T-v_4v_5$
containing $v_5$ and $D'$ be a $\gamma^{d}_2$-set of $T'$ contains
no leaf. Observe that $v_6\in D'$ and $D'\cup \{w_1, w_2, \cdots,
w_t\}$ is a $2DD$-set of $T$. It follows from $n(T)=n(T')+2t+2$,
$l(T)=l(T')+t-1$ and $s(T)=s(T')+t-1$ that $\gamma^{d}_2(T)\leq
\gamma^{d}_2(T')+t\leq
\frac{n(T')+s(T')+l(T')}{4}+t=\frac{n(T)-2t-2+s(T)-t+1+l(T)-t+1}{4}+t=\frac{n(T)+s(T)+l(T)}{4}$.
Suppose that $\gamma^{d}_2(T)=\frac{n(T)+s(T)+l(T)}{4}$, then we
have equality throughout the above inequality chain. In particular,
$\gamma^{d}_2(T')=\frac{n(T')+s(T')+l(T')}{4}$. By the inductive
hypothesis, $(T', S^{*})\in \mathscr{T}_2$ for some labeling
$S^{*}$. Since $v_5$ is a leaf in $T'$, by Observation~4.1(a), it
has status $C$, and then $v_6$ has status $A$. Let $S$ be obtained
from the labeling $S^{*}$ by labeling the vertices $v_3, v_4$ with
label $B, D$, respectively. And moreover, labeling $w_1, w_2,
\cdots, w_t$ with label $A$, and label their leaf-neighbors with
label $C$. Then, $(T, S)$ can be obtained from $(T', S^{*})$ by
doing the operation $\mathscr{O}_3$ for one time and the operation
$\mathscr{O}_2$ for $t-1$ times. Thus, $(T, S)\in \mathscr{T}_2$.

Next we consider the case of $a\geq 1$. Let $u_1, u_2, \cdots, u_a$
be all neighbors of $v_6$ outside $P$ and $u_i'$ be the
leaf-neighbor of $u_i$ ($i=1, 2, \cdots, a$). Let $T'=T-\{u_1, u_2,
\cdots, u_a, u_1', u_2', \cdots, u_a'\}$ and $D'$ be a
$\gamma^{d}_2$-set of $T'$ contains no leaf. Note that $v_6$ has
degree two in $T'$, and $D'\cup \{u_1, u_2, \cdots, u_a\}$ is a
$2DD$-set of $T$. It follows from $n(T)=n(T')+2a$, $l(T)=l(T')+a$
and $s(T)=s(T')+a$ that $\gamma^{d}_2(T)\leq \gamma^{d}_2(T')+a\leq
\frac{n(T')+s(T')+l(T')}{4}+a=\frac{n(T)-2a+s(T)-a+l(T)-a}{4}+a=\frac{n(T)+s(T)+l(T)}{4}$.
Suppose that $\gamma^{d}_2(T)=\frac{n(T)+s(T)+l(T)}{4}$, then we
have equality throughout the above inequality chain. In particular,
$\gamma^{d}_2(T')=\frac{n(T')+s(T')+l(T')}{4}$. By the inductive
hypothesis, $(T', S^{*})\in \mathscr{T}_2$ for some labeling
$S^{*}$.

If $t\geq 2$, by Lemma~4.3, there exists a set $D_1$ with order
$\frac{n(T')+s(T')+l(T')}{4}-1$ such that each vertex of $T'$ is
$2D$-dominated by $D_1$ except for $v_1$, and $v_3$ belongs to
$D_1$. Since leaf-neighbor of each $w_i$ ($i=2, 3, \cdots, t$) is
$2D$-dominated by $D_1$, without loss of generality, we can assume
that each $w_i$ ($i=2, 3, \cdots, t$) belongs to $D_1$. Note that
$d(v_4)=d(v_5)=d(v_6)=2$ in $T'$ and $\{v_4, v_5, v_6, v_7\}\cap
D_1\neq \emptyset$, we construct a set $D_2=(D_1\setminus \{v_4,
v_5, v_6\})\cup \{v_7\}$, each vertex of $T'$ is $2D$-dominated by
$D_2$ except for $v_1$ and $|D_2|\leq |D_1|$. Let $D_3$ be a set
which is obtained from $D_2$ by deleting $v_3$, and adding all
neighbors of $v_6$ outside $P$ and $v_2$. It is easy to see that
$D_3$ is a $2DD$-set of $T$, and $|D_3|\leq
\frac{n(T)+s(T)+l(T)}{4}-1$, it is impossible.

If $t=1$, the vertices $v_1$ and $v_2$ have status $C$ and $A$,
respectively, in $S^{*}$. And so, $v_3$ has status $A$ or $B$.

In the former case, it follows from $d(v_1)=d(v_2)=d(v_3)=d(v_4)=2$
and the definition of $\mathscr{T}_2$ that $v_1v_2v_3v_4$ is the
basic path of $T'$, and then $v_4$ has status $C$. Moreover, $v_5,
v_6$ have status $D, B$, respectively. Let $S$ be obtained from the
labeling $S^{*}$ by labeling each $u_i$ with label $A$, and each
$u_i'$ with label $C$. Then, $(T, S)$ can be obtained from $(T',
S^{*})$ by doing the operation $\mathscr{O}_2$ for $a$ times. Thus,
$(T, S)\in \mathscr{T}_2$.

In the latter case, from the definition of $\mathscr{T}_2$, $v_4,
v_5, v_6$ have status $D, C, A$, respectively. And $v_7$ has status
$A$ or $B$. Assume that sta$(v_7)=A$. If $d(v_7)=2$, we have that
$v_5v_6v_7v_8$ is the basic path of $T'$. Let $S^{*}_1$ be obtained
from $S^{*}$ by changing the status $v_3, v_4, v_5, v_6$ to $A, C,
D, B$, respectively, and clearly, $(T', S^{*}_1)\in \mathscr{T}_2$.
Let $S$ be obtained from the labeling $S^{*}_1$ by labeling each
$u_i$ with label $A$, and each $u_i'$ with label $C$. Then, $(T, S)$
can be obtained from $(T', S^{*}_1)$ by doing the operation
$\mathscr{O}_2$ for $a$ times. Thus, $(T, S)\in \mathscr{T}_2$. If
sta$(v_7)=A$ and $d(v_7)\geq 3$, or sta$(v_7)=B$, let $S$ be
obtained from the labeling $S^{*}$ by labeling each $u_i$ with label
$A$, and each $u_i'$ with label $C$. Then, $(T, S)$ can be obtained
from $(T', S^{*})$ by doing the operation $\mathscr{O}_4$ for $a$
times. Thus, $(T, S)\in \mathscr{T}_2$.
 \ep


\end{document}